\documentstyle[amscd,amssymb,verbatim,12pt]{amsart}
\newcommand{\OO}{{\cal O}}

\newcommand{\MM}{{\cal M}}

\newcommand{\mg}{{\frak m}}

\newcommand{\Proj}{{\Bbb P}}

\newcommand{\Th}{\Theta}

\newcommand{\Pic}{\operatorname{Pic}}

\numberwithin{equation}{subsection}

\newtheorem{thm}{Theorem}[section]
\newtheorem{prop}[thm]{Proposition}
\newtheorem{lem}[thm]{Lemma}
\newtheorem{cor}[thm]{Corollary}

\newcommand{\Pf}{\noindent {\it Proof}}
\newcommand{\id}{\operatorname{id}}

\newcommand{\rk}{\operatorname{rk}}
\newcommand{\ra}{\rightarrow}

\newcommand{\PP}{{\cal P}}

\newcommand{\SS}{{\cal S}}
\newcommand{\LL}{{\cal L}}

\newcommand{\Hom}{\operatorname{Hom}}
\newcommand{\Ext}{\operatorname{Ext}}

\newcommand{\Res}{\operatorname{Res}}

\renewcommand{\a}{\alpha}

\newcommand{\om}{\omega}
\newcommand{\De}{\Delta}

\newcommand{\th}{\theta}
\newcommand{\C}{{\Bbb C}}

\newcommand{\wt}{\widetilde}

\newcommand{\ed}{\qed\vspace{3mm}}

\title{Triple Massey products on curves, Fay's trisecant identity
and tangents to the canonical embedding}
\author{A. Polishchuk}
\thanks{This work is partially supported by the NSF grant
DMS-0070967}
\subjclass{Primary 14H42; Secondary 15A15}

\begin{document}
\begin{abstract}
We show that Fay's trisecant identity follows from the $A_{\infty}$-constraint
satisfied by certain triple Massey products in the derived category of
coherent sheaves on a curve. We also deduce the matrix analogue of this
identity that can be conveniently formulated using quasideterminants
of matrices with non-commuting entries. On the other hand, looking at more
special Massey products we derive a formula for the tangent line to a
canonically embedded curve at a given point. 
\end{abstract}

\maketitle

\bigskip

Fay's trisecant identity is an important
special identity satisfied by theta functions on the
Jacobian of a curve $C$ (see e.g. \cite{A} for an account
of its relation with Schottky problem). Here is a list
of some works containing its proof: \cite{BV}, \cite{BL}, \cite{FK},
\cite{F}, \cite{Gun}, \cite{Poor}, \cite{Raina}, 
(characteristic zero); \cite{Kempf} (arbitrary characteristic).
In this paper we give a new proof of this identity
(valid in arbitrary characteristic).
It turns out that Fay's trisecant identity
follows from the $A_{\infty}$-identity satisfied by certain
triple Massey products on $C$. However, we should stress that our proof does not
dependent on the general theory of Massey products. We
compute explicitly the triple products we need, and the reader can take the
answer as a definition. The main identity between these triple products
is an easy consequence of the residue theorem. The trisecant identity 
follows immediately once we express our
Massey products in terms of theta functions using the Riemann's theorem.
Looking at similar Massey products associated with vector bundles,
we recover the matrix analogue of the trisecant identity
(involving the so called Cauchy-Szeg\"o kernels)
obtained in \cite{BV} and \cite{F2},\cite{F3}.
We observe that this matrix trisecant identity can be
conveniently written using the notion of quasideterminant
introduced by I.~Gelfand and V.~Retakh (see \cite{GR}).
On the other hand, identities with some more special triple Massey
products lead to a relation between tangents
to the canonical embedding of $C$ at triples of points.
As a corollary, we obtain a formula for the tangent
line to a canonically embedded curve at a given point.

In the case when $C$ is an elliptic curve,
the trisecant identity is equivalent to the associative Yang-Baxter
equation satisfied by the Kronecker function (see \cite{P-CYBE}).
In this case the Massey products we consider are the same as in
{\it loc.cit.}. Using the homological mirror symmetry for an elliptic
curve one can express these products in terms of indefinite theta
series (see \cite{P-MP}).
The category of coherent sheaves on a curve $C$ of higher genus
can be considered as a subcategory in the coherent sheaves on
the Jacobian $J$ of $C$. It would be interesting to study implications
for our Massey products
of the homological mirror symmetry for $J$.

\section{A family of triple Massey products}

In this section we look at a family of triple Massey products
involving line bundles and structure sheaves of points on a curve $C$.

\subsection{Definition}

For general definitions concerning triple Massey products in the context
of derived categories we refer to \cite{GM} (Exercises to IV.2) and 
\cite{P-MP}.

In this paper we consider triple Massey products in the derived
category of coherent sheaves on a curve $C$ corresponding to the
following triples of composable morphisms:
$$\OO_C\ra\OO_x\stackrel{[1]}{\ra} L\ra\OO_y,$$
where $\OO_x$ and $\OO_y$ are structure sheaves of points,
$L$ is a line bundle on $C$, the arrow with $[1]$ is a morphism of degree
$1$ (i.e., a morphism $\OO_x\ra L[1]$).

Recall that to compute such a Massey product we have to start by taking
the cone of the middle arrow $\OO_x\ra L[1]$. Assuming that this morphism is
non-zero we get a distinguished triangle
$$\OO_x\ra L[1]\ra L(x)[1]\ra\OO_x[1].$$
Using this triangle we have to lift morphisms $\OO_C\ra\OO_x$ and $L\ra\OO_y$
to morphisms $\OO_C\ra L(x)$ and $L(x)\ra\OO_y$, respectively.
Now the triple Massey product is the composition $\OO_C\ra\OO_y$
of the obtained two morphisms. Note that this operation in general
is ill-defined (liftings not necessarily exist) and multivalued
(liftings are not necessarily unique).

We say that a point $p$ is a base point of a line bundle $M$
if the natural map $H^0(C, M)\ra M|_p$ is zero (e.g., this
is true for any $p$ if $H^0(C,M)=0$).
Let us call a triple $(L,x,y)$
{\it good} if the following three conditions are satisfied:
\begin{enumerate}
\item $x\neq y$;
\item $x$ is a base point of $\om_C L^{-1}$;
\item $y$ is a base point of $L$.
\end{enumerate}
It is easy to see that a triple $(L,x,y)$ is good precisely
when the above Massey product is well-defined and univalued.
In this case using the natural trivializations of $\Hom(\OO_C,\OO_x)$,
$\Hom(\OO_C,\OO_y)$ and the isomorphisms
$\Ext^1(\OO_x,L)\simeq (\om_C^{-1} L)|_x$,
$\Hom(L,\OO_y)\simeq L^{-1}|_y$, we can consider the Massey product
associated with $(L,x,y)$ as an element
$$m_3(L,x,y)\in (\om_C L^{-1})|_x\otimes L|_y.$$
Here is an explicit description of this element.
Since $x$ is a base point of $\om_C L^{-1}$ we have
$h^0(\om_C L^{-1})=h^0(\om_C L^{-1}(-x))$. Therefore,
$h^0(L(x))=h^0(L)+1$. So we can choose an element
$s\in H^0(C,L(x))$ that does not belong to the subspace $H^0(L)$.
In other words, $s$ can be considered as a rational section of
$L$ with a pole of order $1$ at $x$, such that the corresponding
residue $\Res_x(s)\in L(x)|_x\simeq (\om_C^{-1}L)|_x$ is not zero.
On the other hand, we can evaluate $s$ at the point $y$
to get an element $s(y)\in L|_y$. Furthermore, $s(y)$ does not
change if we add to $s$ an element of $H^0(L)$, since $y$ is a base point
of $L$. Hence, the ratio
$$\frac{s(y)}{\Res_x(s)}\in (\om_C L^{-1})|_x\otimes L|_y$$
does not depend on any choices made. It is easy to see that this element
coincides with $m_3(L,x,y)$.

From this description it is clear that $m_3(L,x,y)=0$ if and only if
$y$ is a base point of $L(x)$, or equivalently, $h^0(L(x-y))=h^0(L)+1$.
Note that if the line bundle $L'$ is isomorphic to $L'$
then the corresponding spaces
$(\om_C L^{-1})|_x\otimes L|_y$ and $(\om_C(L')^{-1})|_x\otimes L'|_y$
can be canonically identified, and we have
$m_3(L',x,y)=m_3(L,x,y)$.

In the case when the genus of $C$ is zero, a triple
$(L,x,y)$ is good if and only if $L\simeq\OO(-1)$ and $x\neq y$.
The corresponding Massey product gives a trivialization of
$\OO(-1)\boxtimes\OO(-1)$ outside the diagonal.

In the case $g\ge 1$ we have the theta divisor $\Th\subset\Pic^{g-1}(C)$
(the locus where $h^0$ is non-zero)
and for every line bundle $L$ of degree $g-1$ such that $L\not\in\Th$
the triple $(L,x,y)$ is good whenever $x\neq y$. For such a triple
one has $m_3(L,x,y)=0$ if and only if $L(x-y)\in\Th$.

\subsection{Identity}\label{masseyid}

\begin{thm}\label{idthm} Let $x_1,\ldots,x_n$ be a collection of points,
$L_1,\ldots,L_n$ be a collection of line bundles on $C$ such that
$L_1\ldots L_n\simeq\om_C$. Let
$\a\in H^0(C,\om_C L_1^{-1}\ldots L_n^{-1})$
be a trivialization. Then one has
$$\sum_{i=1}^n \a(x_i)\cdot \prod_{j\neq i} m_3(L_j,x_j,x_i)=0$$
provided that all triples $(L_j,x_j,x_i)$ for $j\neq i$ are good.
\end{thm}

\Pf . For every $j$ choose $s_j\in H^0(C,L_j(x_j))\setminus H^0(C,L_j)$.
Then we can consider $s_1\ldots s_n$ as a rational $1$-form with
poles of order $1$ at $x_1,\ldots, x_n$ (note that all these points are
distinct). Applying the residue theorem we obtain
$$\sum_{i=1}^n \Res_{x_i}(s_i)\prod_{j\neq i}s_j(x_i)=0.$$
Dividing by $\prod_i \Res_{x_i}(s_i)$ we get the result.
\ed

In the case $n=2$ the above theorem simply tells that
if $L_1L_2\simeq\om_C$ then
\begin{equation}\label{skewsym}
\a(y)\cdot m_3(L_1,x,y)=-\a(x)\cdot m_3(L_2,y,x),
\end{equation}
where $\a\in H^0(C,\om_C L_1^{-1}L_2^{-1})$ is the corresponding
trivialization.
The case $n=3$ is equivalent to the $A_{\infty}$-axiom applied to
the sequence of $5$ composable morphisms
$$\OO_C\ra\OO_x\stackrel{[1]}{\ra} L\ra\OO_y\stackrel{[1]}{\ra} M\ra\OO_z.$$
Indeed, the compatibility of the triple products with Serre duality
(see \cite{P-hmc}) implies that the Massey product 
$$\OO_x\stackrel{[1]}\ra L\ra\OO_y\stackrel{[1]}\ra M$$
essentially coincides with $m_3(ML^{-1},y,x)$. Now using the
skew-symmetry (\ref{skewsym}), one can easily see that the
above $A_{\infty}$-axiom is equivalent to the identity of the Theorem
\ref{idthm} applied to $x_1=x$, $x_2=y$, $x_3=z$, $L_1=L$, $L_2=ML^{-1}$,
$L_3=\om_C M^{-1}$.


\begin{cor}\label{maincor} Assume that $g\ge 1$.
Let $L_1$ and $L_2$ be line bundles of degree $g-1$ with
$H^0(C,L_1)=H^0(C,L_2)=0$, $z_1,\ldots,z_n$,
$t_1,\ldots,t_n$ be points on $C$ such that
$L_1L_2\simeq\om_C(\sum_{i=1}^n z_i-\sum_{i=1}^n t_i)$.
Let us fix a trivialization $\a\in H^0(C,\om_C L_1^{-1}L_2^{-1}
(\sum_i(z_i-t_i)))$.
Then for a pair of points $x,y\in C$ one has
$$\sum_{i=1}^n\Res_{z_i}(\a)\cdot m_3(L_1,x,z_i)m_3(L_2,y,z_i)+
\a(x)\cdot m_3(L_2,y,x)+\a(y)\cdot m_3(L_1,x,y)=0.$$
provided that the points $(x,y,z_1,\ldots,z_n)$ are distinct and
disjoint from the set $\{t_1,\ldots,t_n\}$.
\end{cor}

\Pf . Apply Theorem \ref{idthm} to the line bundles
$L_1$, $L_2$ and $L_{2+i}=\OO_C(t_i-z_i)$ and the points $x,y,z_1,\ldots,z_n$,
and then use the equality
$m_3(\OO_C(t-z),z,y)=1$ for $y\neq t$.
\ed

Here is a version of the case $n=1$ of
the above Corollary for line
bundles $L_1$ and $L_2$ with $h^0=1$.

\begin{cor}\label{cor2} Assume that $g\ge 3$.
Let $E_1,E_2$ be effective divisors of degree $g-3$, $x,y,z,t$
be distinct points disjoint from $E_1$ and $E_2$, such that
$\OO_C(E_1+E_2+x+y+z+t)\simeq\om_C$. Assume also that
$h^0(E_1+y+z)=h^0(E_2+x+z)=1$.
Let $\eta$ be a non-zero
$1$-form on $C$ with the zero divisor $E_1+E_2+x+y+z+t$. Then
\begin{align*}
&\eta'(z)m_3(\OO_C(E_1+y+z),x,z)m_3(\OO_C(E_2+x+z),y,z)+\\
&\eta'(x)m_3(\OO_C(E_2+x+z),y,x)+\eta'(y)m_3(\OO_C(E_1+y+z),x,y)=0.
\end{align*}
Here $\eta'(x)$ denotes the reduction of $\eta\in\mg_x\om_C$
modulo $\mg_x^2\om_C$, where $\mg_x$ is a maximal ideal of $x$;
it is an element of $\mg_x/\mg_x^2\otimes\om_C=\om_C^2|_x$.
\end{cor}

\Pf . First we note that due to the assumption on $h^0$, all
the triple products involved are well-defined and univalued.
It remains to apply Theorem \ref{idthm} to line bundles
$L_1=\OO_C(E_1+y+z)$, $L_2=\OO_C(E_2+x+z)$, $L_3=\OO_C(t-z)$ and points
$x,y,z$.
\ed

\section{Theta functions}

In this section we assume that $g\ge 1$. First we recall
Riemann's theorem in the form we need and then we express
the Massey products from the previous section in terms of theta
functions. The identity between Massey products
considered above translates into Fay's trisecant identity.

\subsection{Riemann's theorem}

Let us choose a base point $p_0\in C$.
Let $J$ be the Jacobian of $C$, $i:C\ra J$ be the embedding sending
a point $p$ to $\OO_C(p-p_0)$. We will identify $C$ with a subvariety
of $J$ by means of this embedding. For a pair of points $x,y\in C$
the difference $i(x)-i(y)$ will be denoted simply by $x-y$ (it
does not depend on a choice of $p_0$).

For every line bundle $L$ of degree $g-1$ on $C$ 
there is an associated theta-divisor $\Th_L\subset J$. As a set
this is the locus of $\xi\in J$ such that
$h^0(L\otimes\xi)\neq 0$.
We have $\Th_{L\otimes\xi}=\Th_L-\xi$ for every $\xi\in J$,
hence, all the divisors $\Th_L$ are translations of
each other.
For a point $\xi\in J$ we define the line bundle $\PP_{\xi}$
on $J$ by
$$\PP_{\xi}=t_{\xi}^*\LL\otimes \LL^{-1}\otimes\LL^{-1}|_{\xi}
\otimes\LL|_0$$
where $\LL=\OO_J(\Th_L)$, $t_{\xi}:J\ra J$ is the
translation by $\xi$. The theorem of the cube implies
that the right-hand side
does not depend on a choice of theta-divisor $\Th_L$.

The following theorem is a reformulation of
the classical theorem due to Riemann. Apparently, this
reformulation is well-known but we could not locate the
exact reference.

\begin{thm}\label{Riemann} 
One has the following isomorphism of line bundles on $C$:
$$\OO_J(\Th_L)|_C\simeq\om_C L^{-1}(p_0).$$
\end{thm}

\Pf . Recall that there is a natural embedding of 
$C$ into $\hat{J}$ given by the universal family on
$C\times J$ trivialized over $p_0\times J$ (in fact, $\hat{J}\simeq J$
but we don't use this isomorphism here to avoid confusion).
Thus, we can consider bundles on $C$ as coherent sheaves on $\hat{J}$
and apply the Fourier-Mukai transform $\SS$ to them to get objects
of the derived category of coherent sheaves on $J$ (see \cite{M}).
It is easy to see that $\OO_J(\Th_L)\simeq\det^{-1}(\SS(L))$.
On the other hand, by the base change formula we get
$$\SS(L)|_C\simeq R\pi_{2*}(\pi_1^*(L(-p_0))(\De_C))(-p_0),$$
where $\De_C$ is the image of the diagonal
embedding $\De:C\ra C\times C$,
$\pi_1,\pi_2:C\times C\ra C$ are projections. 
Let us denote $L_1=L(-p_0)$. Note that $\rk\SS(L)=\chi(L)=0,$ 
so we have
$$\det\SS(L)|_C\simeq \det(R\pi_{2*}(\pi_1^*L_1(\De_C))).$$
Now from the exact sequence of sheaves on $C\times C$
$$0\ra \pi_1^*L_1\ra \pi_1^*L_1(\De_C)\ra
\De_*(L_1\otimes \om_C^{-1})\ra 0$$
we get
$$\det R\pi_{2*}(\pi_1^*L_1(\De_C))\simeq L_1\otimes \om_C^{-1}$$
as required.
\ed

The classical corollary of the above theorem is that for fixed
$L\in\Pic^{g-1}(C)$ the map
$\xi\mapsto (\Th_{L}+\xi)\cap C=\Th_{L\otimes\xi^{-1}}\cap C$
is a birational isomorphism from $J$ to $S^g(C)$ which is inverse
(up to translation) to the natural map
$S^g(C)\ra\Pic^g(C):D\mapsto\OO_C(D)$.

The above theorem also immediately implies that for every $\xi\in J$
the restriction of $\PP_{-\xi}$ to $C$ is isomorphic to $\xi$
(considered as a line bundle on $C$).
Let us denote by $\th_L$ a global section of a line bundle on $J$ 
which has $\Th_L$ as the zero locus (it is defined uniquely up to
a non-zero constant). 
We will call functions of the form $\th_L$ where $\deg(L)=g-1$,
{\it theta functions of degree $1$ on $J$}.
We conclude that if $\th_L|_C$ does not vanish identically, then
$$\frac{t_{-\xi}^*\th_L}{\th_L}|_{C}$$
is a rational section of a line bundle isomorphic to
$\xi$ with the divisor of poles $\Th_L\cap C$.

\subsection{Expression of Massey products in terms of theta functions}
\label{exprthetasec}

Let $D_1$ and $D_2$ be a pair of effective divisors on $C$ such that
$\om_C\simeq\OO_C(D_1+D_2)$. We denote by $\eta$ a non-zero global
$1$-form on $C$
with zeroes at $D_1+D_2$ (it is defined uniquely up to a non-zero
constant). Let us assume that $h^0(D_1)=h^0(D_2)=1$.     
We are going to calculate the Massey products
$m_3(\xi(D_2),x,y)$ for generic $\xi\in J$, $x,y\in C$
in terms of the theta function associated with $D_1$.
Note that since $x$ and $y$ are not contained in the support of $D_2$,
we can consider $m_3(\xi(D_2),x,y)$ as an element of
$\xi|_y\otimes\xi^{-1}|_x\otimes\om_C|_x$. 
To compute it we have to construct a rational section of $\xi$
with $D_2+x$ as the divisor of poles. As we have seen above,
for this we need to represent $D_2+x$ in the form $\Th_L\cap C$
for some line bundle $L$ of degree $g-1$. Let us set $p_0=x$,
so that the embedding of $C$ into $J$ is given by $p\mapsto\OO_C(p-x)$.
Then we claim that $L=\OO_C(D_1)$ will work.
Indeed, by Theorem \ref{Riemann} we have
$$\OO_J(\Th_{D_1})|_C\simeq\om_C(x-D_1)\simeq\OO_C(D_2+x).$$
Furthermore, for generic $x$ and $y$ we have $h^0(D_1+y-x)=0$
and $h^0(D_2+x)=1$. The former condition implies that $\Th_{D_1}$
does not contain $y-x$, while the latter condition implies
that $\Th_{D_1}\cap C$ is precisely $D_2+x$.
Thus,
$$s=\frac{t^*_{-\xi}\th_{D_1}}{\th_{D_1}}|_C$$
is a rational section of a line bundle isomorphic to $\xi$ with
$D_2+x$ as the divisor of poles. By abuse of notation we consider
$s$ as a rational section of $\xi$ itself
(this does not lead to ambiguity 
since the expression for $m_3$ is invariant
under rescaling of $s$).
Now we have
$$
m_3(\xi(D_2),x,y)=\frac{s(y)}{\Res_x(s)}=
\frac{\th_{D_1}(y-x-\xi)}
{\th_{D_1}(y-x)\th_{D_1}(-\xi)
\Res_{x}\frac{1}{t^*_{-i(x)}\th_{D_1}|_C}}
$$
Using the natural identification of the cotangent space to $J$
at $0$ with $H^0(C,\om_C)$ we can consider the derivative
$\th'_{D_1}(0)$ as a global $1$-form on $C$ (more precisely,
it is a $1$-form with values in a one-dimensional vector space,
the stalk of $\OO_J(\Th_{D_1})$ at $0$). Then it is easy to see
that 
$$\Res_{x}\frac{1}{t^*_{-i(x)}\th_{D_1}|_C}=\frac{1}{\th'_{D_1}(0)(x)}.$$
Therefore, we arrive to the following expression for the Massey product.

\begin{lem} Under the canonical identification of spaces
$$(\om_C\otimes\xi^{-1})|_x\otimes\xi|_y\simeq \PP_{-\xi}|_{y-x}\otimes
\om_C|_x$$
one has
\begin{equation}\label{masseytheta}
m_3(\xi(D_2),x,y)=
\frac{\th_{D_1}(y-x-\xi)\th'_{D_1}(0)(x)}
{\th_{D_1}(y-x)\th_{D_1}(-\xi)}
\end{equation}
\end{lem}

Now we are going to substitute the above expression for $m_3$ into
the identities we derived in section \ref{masseyid}.
First, we have the skew-symmetry (\ref{skewsym}):
$$\eta(y)\cdot m_3(\xi(D_2),x,y)=
-\eta(x)\cdot m_3(\xi^{-1}(D_1),y,x).$$
Since $\Th_{D_1}=-\Th_{D_2}$, we can normalize $\th_{D_1}$
and $\th_{D_2}$ in such a way that $\th_{D_1}(-\xi)=\th_{D_2}(\xi)$
(when we consider $\xi$ as an element of $J$ we use the
additive notation $-\xi$ for the opposite element).
Then the above identity reduces to
$$\eta(y)\th'_{D_1}(0)(x)=\eta(x)\th'_{D_1}(0)(y).$$
In other words, the global $1$-forms $\eta$ and $\th'_{D_1}(0)$
are proportional. Thus, we recover the well-known fact that
the $1$-form $\th'_{D_1}(0)$ vanishes on $D_1$
(see \cite{F}, ch. I, Cor. 1.4).

Now we are going to look at the identity obtained from Corollary
\ref{maincor}.

\begin{thm}\label{trisecthm}
Let $F(\xi_1,\xi_2)=\frac{\th(\xi_1+\xi_2)}
{\th(\xi_1)\th(\xi_2)}$ where $\xi_1,\xi_2\in J$, $\th$
is a theta-function on $J$ of degree $1$ such that $\th(0)=0$.
Then one has
\begin{equation}\label{mainid}
\begin{align*}
&
\sum_{i=1}^n
(\prod_{j\neq i}F(z_i-z_j,z_j-t_j))\cdot
F(z_i-x,\xi)F(y-z_i,\sum_k (z_k-t_k)+\xi)+\\
&(\prod_{i=1}^n F(x-z_i,z_i-t_i))\cdot F(y-x,\sum_k (z_k-t_k)+\xi)-\\
&(\prod_{i=1}^n F(y-z_i,z_i-t_i))\cdot F(y-x,\xi)=0.
\end{align*}
\end{equation}
where $x,y,z_1,\ldots,z_n,t_1,\ldots,t_n\in C$, $\xi\in J$.
\end{thm}

\Pf . Applying Corollary \ref{maincor} to $L_1=\PP_{\xi}|_C(D_1)$,
$L_2=\PP_{-\xi+\sum_i(t_i-z_i)}|_C(D_2)$ and $\a=\phi\cdot\eta$, where
$\phi\in H^0(C,\PP_{\sum_i(z_i-t_i)}|_C(\sum_i(z_i-t_i)))$
is the corresponding trivialization, we get
\begin{align*}
&\sum_i\Res_{z_i}(\phi\eta)\cdot
\frac{\th_{D_2}(z_i-x+\xi)\th'_{D_2}(0)(x)}
{\th_{D_2}(z_i-x)\th_{D_2}(\xi)}\cdot
\frac{\th_{D_1}(z_i-y+\sum_k(t_k-z_k)-\xi)\th'_{D_1}(0)(y)}
{\th_{D_1}(z_i-y)\th_{D_1}(\sum_k(t_k-z_k)-\xi)}+\\
&\phi(x)\eta(x)\cdot\frac{\th_{D_1}(x-y+\sum_k(t_k-z_k)-\xi)\th'_{D_1}(0)(y)}
{\th_{D_1}(x-y)\th_{D_1}(\sum_k(t_k-z_k)-\xi)}+
\phi(y)\eta(y)\cdot\frac{\th_{D_2}(y-x+\xi)\th'_{D_2}(0)(x)}
{\th_{D_2}(y-x)\th_{D_2}(\xi)}=0.
\end{align*}
As we have seen above, we can take $\eta=\th'_{D_2}(0)=-\th'_{D_1}(0)$.
Now using the equality $\th_{D_1}=[-1]^*\th_{D_2}$ we get
\begin{align*}
&\sum_i\Res_{z_i}(\phi\eta)\cdot
\frac{\th_{D_2}(z_i-x+\xi)}
{\th_{D_2}(z_i-x)\th_{D_2}(\xi)}\cdot
\frac{\th_{D_2}(y-z_i+\sum_k(z_k-t_k)+\xi)}
{\th_{D_2}(y-z_i)\th_{D_2}(\sum_k(z_k-t_k)+\xi)}+\\
&\phi(x)\cdot
\frac{\th_{D_2}(y-x+\sum_k(z_k-t_k)+\xi)}
{\th_{D_2}(y-x)\th_{D_2}(\sum_k(z_k-t_k)+\xi)}-
\phi(y)\cdot
\frac{\th_{D_2}(y-x+\xi)}
{\th_{D_2}(y-x)\th_{D_2}(\xi)}=0.
\end{align*}
We can choose $\phi$ in the form $\phi=\prod_j\phi_{t_j,z_j}$
where $\phi_{t_i,z_i}$ is a rational section of $\PP_{z_i-t_i}|_C$
with pole of order $1$ at $z_i$ and zero of order $1$ at $t_i$.
Assuming that $t_i$ and $z_i$ are generic, we can take
$$\phi_{t_i,z_i}(t)=
\frac{\th_{D_2}(t-t_i)}{\th_{D_2}(t-z_i)}.$$
Indeed, clearly this is a section of
a line bundle isomorphic to $\PP_{z_i-t_i}|_C$.
On the other hand, for generic $z_i,t_i$ we have
$\th_{D_2}(t_i-z_i)\neq 0$, while
$D_1+z_i$ is the only effective divisor in its
linear series. Therefore, the denominator of $\phi_{t_i,z_i}$
vanishes presicely on $D_1+z_i$. Its numerator vanishes on
$D_1+t_i$, hence $\phi_{t_i,z_i}$ has the only pole of order $1$ at $z_i$.
Now we have
$$\Res_{z_i}(\prod_j\phi_{t_j,z_j}\eta)=
\prod_{j\neq i}\frac{\th_{D_2}(z_i-t_j)}{\th_{D_2}(z_i-z_j)}\cdot
\th_{D_2}(z_i-t_i).$$
Thus, we can rewrite the above identity as follows (where $\th=\th_{D_2}$):
\begin{align*}
&\sum_i
\prod_{j\neq i}\frac{\th(z_i-t_j)}{\th(z_i-z_j)}\cdot
\th(z_i-t_i)\cdot F(z_i-x,\xi)F(y-z_i,\sum_k(z_k-t_k)+\xi)+\\
&\prod_i\frac{\th(x-t_i)}{\th(x-z_i)}\cdot F(y-x,\sum_k(z_k-t_k)+\xi)
-\prod_i\frac{\th(y-t_i)}{\th(y-z_i)}\cdot F(y-x,\xi)=0.
\end{align*}
Dividing by $\prod_i\th(z_i-t_i)$ we get the required identity for
$\th=\th_{D_2}$. Our only assumption on $D_2$ was that $h^0(D_2)=1$
which is an open condition. Thus, the same identity holds if
$D_2$ is an arbitrary effective divisor of degree $g-1$, which means
that $\th$ can be any theta function of degree $1$ such that $\th(0)=0$.
\ed

In the case of elliptic curves the function $F$ coincides with the Kronecker
function studied in \cite{Kr},\cite{Weil},\cite{Z}. The equation
(\ref{mainid})
in this case reduces to the scalar associative Yang-Baxter
equation studied in \cite{P-CYBE}.
Thus, the functions $F(y-x,\xi)$ on (a covering of) $C\times C\times J$
should be considered as the higher genus analogues of the Kronecker
function.

The identity (\ref{mainid}) for $n=1$ can be rewritten in the following form:
\begin{align*}
&\frac{\th(x-t)\th(y-z)}{\th(x-z)\th(y-t)}\cdot\th(\xi)\th(\xi+y-x+z-t)+\\
&\frac{\th(z-t)\th(y-x)}{\th(z-x)\th(y-t)}\cdot\th(\xi+z-x)\th(\xi+y-t)=\\
&\th(\xi+z-t)\th(\xi+y-x)
\end{align*}
where $\th$ is a theta function of degree $1$ vanishing at $0$.
This identity is equivalent to the trisecant identity (see \cite{FK},VII.6).
In fact, the identity (\ref{mainid}) for general $n$ is a
formal consequence of the case $n=1$, however, the derivation is
not straightforward (see section \ref{bunidsec}).

The following particular case of (\ref{mainid})  deserves attention
because of its symmetric form.

\begin{cor}\label{divisorid}
In the notation of the theorem one has
$$\sum_{i=0}^n\prod_{j\neq i}F(z_i-z_j,z_j-t_j)\cdot F(y-z_i,
\sum_{k=0}^n z_k-\sum_{k=0}^n t_k)=\prod_{i=0}^n F(y-z_i,z_i-t_i),
$$
where $y,z_0,\ldots,z_n,t_0,\ldots,t_n\in C$.
\end{cor}

\Pf . Set $x=z_0$, $\xi=z_0-t_0$ in (\ref{mainid}).
\ed

\subsection{Expression in terms of the prime form}
\label{primeform}

In this section we work over $\C$.
We are going to derive a different formula for the Massey product
$m_3(L,x,y)$ where $\deg(L)=g-1$ involving the Schottky-Klein prime
form. Let us recall its definition. Let us pick a non-singular odd
theta-characteristic on $C$, i.e. an effective divisor $D$ with
$\OO_C(2D)\simeq\om_C$ and $h^0(D)=1$
(it is known that such $D$ always exist, see \cite{F}).
As we have seen above, the non-zero global $1$-form $\eta$ with $2D$
as zero divisor is proportional to $\th'_{D}(0)$. Hence, we can
choose a global section $h_{D}$ of $\OO_C(D)$ satisfying
$$h_{D}^2=\th'_{D}(0).$$
Now the prime-form is defined by the formula
$$E(x,y)=\frac{\th_{D}(y-x)}{h_{D}(x)h_{D}(y)}.$$
Since $\th_{D}$ is odd, we have $E(y,x)=-E(x,y)$.
One can deduce from Riemann's theorem that $E(x,y)$ does not depend on
$D$ and that for any line bundle $L$ of degree $g-1$ 
$$t\mapsto \frac{\th_L(t-x)}{E(x,t)},$$
is a meromorphic section of $L$ with the unique simple pole
at $t=x$ (it is holomorphic if $h^0(L)>0$).
This implies the following formula for the Massey product
$m_3(L,x,y)$.

\begin{lem} Assume that $L$ is a line bundle of degree $g-1$
with $h^0(L)=0$. Then for every pair
of distinct points $x,y\in C$ one has
$$m_3(L,x,y)=\frac{\th_L(y-x)}{E(x,y)\th_L(0)}.$$
\end{lem}

\Pf . By the definition we have
$$m_3(L,x,y)=\frac{s(y)}{\Res_x s},$$
where $s$ is a meromorphic section of $L$ having the unique
simple pole at $x$ with non-zero residue. Since $\th_L(0)\neq 0$
we can take
$s(t)=\frac{\th_L(t-x)}{E(x,t)}$. Thus,
$$m_3(L,x,y)=\frac{\th_L(y-x)}{E(x,y)\th_L(0)
\Res_{t=x}\frac{1}{E(x,t)}}.$$
It remains to notice that
$$\Res_{t=x}\frac{1}{E(x,t)}=\Res_{t=x}\frac{h_{D}(x)h_{D}(t)}
{\th_{D}(t-x)}=\Res_{t=x}\frac{h_{D}(t)^2}{\th_{D}(t-x)}=1.$$
\ed

Using this computation we can rewrite Corollary \ref{maincor} in
the following form.

\begin{thm} Let $\th$ be arbitrary theta function of degree $1$ on $J$.
Then for any collection of generic points
$x$,$y$,$z_1,\ldots,z_n$,$t_1,\ldots,t_n$ in $C$ one has
\begin{equation}\label{another}
\begin{array}{l}
\sum_{i=1}^n\prod_{j\neq i}\frac{E(t_j,z_i)}{E(z_j,z_i)}
\cdot\frac{E(t_i,z_i)E(x,y)}{E(x,z_i)E(y,z_i)}\cdot
\th(z_i-x)\th(y-z_i+\sum_i(z_i-t_i))+\\
\prod_i\frac{E(t_i,y)}{E(z_i,y)}\cdot
\th(y-x)\th(\sum_i(z_i-t_i))=
\prod_i\frac{E(t_i,x)}{E(z_i,x)}\cdot
\th(y-x+\sum_i(z_i-t_i))\th(0).
\end{array}
\end{equation}
\end{thm}

\Pf . We can assume that $\th=\th_L$ where $h^0(L)=0$.
Then we can apply Corollary \ref{maincor} to $L_1=L$,
$L_2=\om_C\otimes L^{-1}(\sum_i(z_i-t_i))$ and
$$\a=\prod_i\frac{E(t_i,\cdot)}{E(z_i,\cdot)}.$$
It remains to use the equalities
$\th_{\om_C L^{-1}}=[-1]^*\th_L$ and
$\Res_{z_i}\frac{1}{E(z_i,\cdot)}=1$.
\ed

In the particular case $n=1$ the equation (\ref{another}) becomes
$$E(x,y)E(t,z)\th(z-x)\th(y-t)
+E(x,z)E(y,t)\th(y-x)\th(z-t)=\\
E(x,t)E(y,z)\th(y-x+z-t)\th(0)$$
where $\th$ is any theta function of degree $1$ on $J$.
This identity is equivalent to another form of the trisecant
identity (the formula (45) in \cite{F}, ch.II).

\section{Vector bundles}

\subsection{Triple Massey products and Cauchy-Szeg\"o kernels}

The generalization of the notion of a good triple to the case of vector
bundles is straightforward. Namely, we say that a triple $(V,x,y)$, where
$V$ is a vector bundle on $C$, is {\it good} if the 
following three conditions are satisfied:
\begin{enumerate}
\item $x\neq y$;
\item $x$ is a base point of $\om_C\otimes V^{\vee}$;
\item $y$ is a base point of $V$.
\end{enumerate}
Here we say that a point $z$ is a base point of a vector bundle $W$ if
the evaluation map $H^0(C,W)\ra W|_z$ is zero.
It is easy to see that for a pair of vector bundles $V_1$ and $V_2$,
the triple $(V_1^{\vee}\otimes V_2,x,y)$ is good precisely
when the triple Massey products of the type 
\begin{equation}\label{bundletype}
V_1\ra\OO_x\stackrel{[1]}\ra V_2\ra\OO_y
\end{equation}
are well-defined and univalued.

For a good triple $(V,x,y)$ an analogue of the triple Massey product
can be defined as follows. Since $x$ is a base point of 
$\om_C\otimes V^{\vee}$, by Serre duality we obtain
that the map $p_x:H^0(C,V(x))\ra V(x)|_x$ is surjective
(with the kernel $H^0(C,V)$).
On the other hand, the natural evaluation map
$H^0(C,V(x))\ra V(x)|_y\simeq V|_y$ vanishes on the kernel of $p_x$.
Therefore, it factors through a unique map
$m_3(V,x,y):V(x)|_x\ra V|_y$. Thus, we can consider $m_3(V,x,y)$ as
a map $V|_x\ra \om_C|_x\otimes V|_y$.
If $V=V_1^{\vee}\otimes V_2$ then this element represents the Massey
product of the type (\ref{bundletype}). Note that every
map $m_3(V,x,y)$ appears as such Massey product (take $V_1=\OO_C$).

Now assume that $\mu(V)=\frac{\deg V}{\rk V}=g-1$ and $h^0(V)=0$
(this implies that $V$ is semistable).
Then for every pair of distinct points $(x,y)$ the triple $(V,x,y)$
is good. In this situation,
the construction of $m_3(V,x,y)$ can be globalized to define a morphism
$$\pi_1^*V\ra \pi_1^*\om_C\otimes \pi_2^* V(\Delta_C)$$
of vector bundles on $C\times C$ (where $\pi_1,\pi_2:C\times C\ra C$ are
projections, $\Delta_C\subset C\times C$ is the diagonal)
whose residue on the diagonal is equal to the identity.
Let $L$ be a line bundle of degree $g-1$ on $C$, and
$\chi$ be a flat vector bundle of rank $r$ on $C$ such that
$h^0(\chi\otimes L)=0$.
Then $m_3(\chi\otimes L,x,y)$ can be considered as a rational section of
a vector bundle on $\MM_r\times C\times C$, where $\MM_r$
is the moduli space of flat bundles of rank $r$ on $C$.
In the particular case when $L$ is a square root of the canonical bundle,
$m_3(\chi\otimes L,x,y)$ coincides
with the Cauchy kernel defined in \cite{BV} and
with the Szeg\"o kernel considered in \cite{F2},\cite{F3} (more precisely,
it corresponds to $S(y,x,\chi)$ in Fay's notation). For
this reason we propose to call it Cauchy-Szeg\"o kernel.

Note that by the definition, the map $m_3(V,x,y)$ is an isomorphism
if and only if $h^0(C,V(x-y))=0$.

\subsection{Identities}\label{bunidsec}

Theorem \ref{idthm} has the following generalization
to the case of vector bundles.

\begin{thm}\label{idthm2} Let $x_1,\ldots,x_n$ be a collection of points,
$V_1,\ldots,V_n$ be a collection of vector bundles on $C$, and
$\a:V_1\otimes\ldots \otimes V_n\ra\om_C$ be a morphism. Then one has
$$\sum_{i=1}^n \a(x_i)(\id_{V_i|_{x_i}}\otimes
(\otimes_{j\neq i} m_3(V_j,x_j,x_i)))=0$$
provided that all triples $(V_j,x_j,x_i)$ for $j\neq i$ are good.
\end{thm}

The proof is a straightforward generalization of the proof of Theorem
\ref{idthm}.


The case $n=2$ of the above theorem leads to the following equality:
\begin{equation}\label{skewsym2}
m_3(V,x,y)=-m_3(V^{\vee}\otimes\om_C,y,x)^*.
\end{equation}
On the other hand, we have the following analogue of
Corollary \ref{maincor}.

\begin{cor}\label{mainbuncor}
Let $V$ be a vector bundle with $\mu(V)=g-1$ such that $h^0(V)=0$,
and let $z_1,\ldots,z_n$, $t_1,\ldots,t_n$ be points on $C$ such that
$h^0(V(\sum_{i=1}^n t_i-\sum_{i=1}^n z_i))=0$. Then for
a pair of points $x,y\in C$ one has
$$\sum_{i=1}^n
m_3(V(\sum_{i=1}^n t_i-\sum_{i=1}^n z_i),z_i,y)\circ m_3(V,x,z_i)+
m_3(V(\sum_{i=1}^n t_i-\sum_{i=1}^n z_i),x,y)-m_3(V,x,y)=0,$$
provided that the points $(x,y,z_1,\ldots,z_n)$ are distinct and
disjoint from the set $\{t_1,\ldots,t_n\}$.
\end{cor}

\Pf . Applying Theorem \ref{idthm2} to the vector bundles
$V_1=V$, $V_2=V^{\vee}\otimes\om_C(\sum_{i=1}^n z_i-\sum_{i=1}^n t_i)$,
$V_{2+i}=\OO_C(t_i-z_i)$, $i=1,\ldots,n$ and the points
$x_1=x,x_2=y,x_{2+i}=z_i$, $i=1,\ldots,n$,
and using the equality
$m_3(\OO_C(t-z),z,y)=1$ for $y\neq t$, we get
\begin{align*}
&\sum_i m_3(V^{\vee}\otimes
\om_C(\sum_{i=1}^n z_i-\sum_{i=1}^n t_i),y,z_i)^*
\circ m_3(V,x,z_i)+\\
&m_3(V^{\vee}\otimes\om_C(\sum_{i=1}^n z_i-\sum_{i=1}^n t_i),y,x)^* +
m_3(V,x,y)=0.
\end{align*}
It remains to use the skew-symmetry (\ref{skewsym2}).
\ed

We will also need the following degenerate version of the case $n=1$ of
the above corollary.

\begin{cor}\label{idcor}
Let $V$ be a vector bundle with $\mu(V)=g-1$ such that $h^0(V)=0$,
and let $x,y,z$ be a triple of distinct points such that
$h^0(V(x-z))=0$. Then
$$m_3(V(x-z),z,y)=m_3(V,x,y)\circ m_3(V,x,z)^{-1}.$$
\end{cor}

\Pf . We apply Theorem \ref{idthm2} to $V_1=V$, $V_2=V^{\vee}\otimes\om_C(z-x)$,
$V_3=\OO_C(x-z)$, $x_1=x$, $x_2=y$, $x_3=z$ and then use the
equalities
$m_3(\OO_C(x-z),z,y)=1$, $m_3(\OO_C(x-z),z,x)=0$, and the
skew-symmetry (\ref{skewsym2}) to get the equality
$$m_3(V(x-z),z,y)\circ m_3(V,x,z)=m_3(V,x,y).$$
Next we note that if $h^0(V(x-z))=0$ then
$m_3(V,x,z)$ is invertible.
\ed

Below we will use the notion of quasideterminant of a matrix
with non-commuting entries, introduced by
I.~Gelfand and V.~Retakh (see \cite{GR} and the references therein).
This notion appears naturally in the problem of inverting matrices
with non-commuting entries. In our case the entries will be linear
maps between vector spaces. 

\begin{thm}\label{quasidetthm}
Let $V$ be a vector bundle with $\mu(V)=g-1$ such
that $h^0(V)=0$, and let $\{x_0\},\{x_1,\ldots,x_n\}$,
$\{y_0\}$ and $\{y_1,\ldots,y_n\}$ be disjoint sets of points such that
$h^0(V(\sum_{i=1}^n x_i-\sum_{i=1}^n y_i))=0$.
Let us consider the $(n+1)\times (n+1)$-matrix $A=(a_{ij})$ with
the entries $a_{ij}=m_3(V,x_j,y_i)$, $0\le i,j\le n$. Then one has
$$m_3(V(\sum_{i=1}^n x_i-\sum_{i=1}^n y_i),x_0,y_0)=|A|_{00}.$$
\end{thm}

\Pf . The case $n=1$ of Corollary \ref{mainbuncor} gives
$$m_3(V(x_1-y_1),x_0,y_0)=m_3(V,x_0,y_0)-m_3(V(x_1-y_1),y_1,y_0)\circ
m_3(V,x_0,y_1).$$
Substituting the expression for $m_3(V(x_1-y_1),y_1,y_0)$ given
by Corollary \ref{idcor}, we get
\begin{equation}\label{quasidetid}
m_3(V(x_1-y_1),x_0,y_0)=m_3(V,x_0,y_0)-m_3(V,x_1,y_0)\circ
m_3(V,x_1,y_1)^{-1}\circ m_3(V,x_0,y_1).
\end{equation}
The expression in the right-hand side is the $(0,0)$-quasideterminant
of the matrix
$$\left(\matrix m_3(V,x_0,y_0) & m_3(V,x_1,y_0) \\
m_3(V,x_0,y_1) & m_3(V,x_1,y_1)\endmatrix\right),$$
so we get the assertion in the case $n=1$.
It remains to iterate this equation and to use the quasideterminant
analogue of the Sylvester identity (see \cite{GR}, Thm. 1.3.1).
\ed

Note that in the proof of the theorem above we only used the case $n=1$
of Corollary \ref{mainbuncor}. Applying this theorem, we see that
the case of arbitrary $n$ is equivalent
to the following identity:
$$|A|_{00}=a_{00}-\sum_{i\ge 1}|A^{i0}|_{0i}\cdot |A^{00}|_{ii}^{-1}
\cdot a_{i0},$$
where $A=(a_{ij})$ is the matrix from the Theorem \ref{quasidetthm},
$A^{ij}$ is the matrix obtained from $A$ by deleting $i$-th row
and $j$-th column. This identity follows easily from the column
expansion and the homological relations for quasideterminants
(see \cite{GR}, Prop. 1.2.7 and Thm. 1.2.3).

In conclusion let us compare the above identity with the
identity obtained by Fay in \cite{F3}.
We fix a line bundle $L$ of degree $g-1$ and set
$$F(\chi,x,y)=m_3(\chi\otimes L,x,y)$$
where $\chi\in\MM_r$ is a flat bundle on $C$, $x,y\in C$. In order
to rewrite the identity (\ref{quasidetid}) as an equality over
$\MM_r\times C^4$ we have to take into account the trivializations
of the flat line bundle $\xi=\OO_C(x_1-y_1)$ at the points $x_0$ and $y_0$.
To construct such trivializations we can use the rational section 
$\phi_{x_1,y_1}=t^*_{-i(x_1)}\th_{D}/t^*_{-i(y_1)}\th_{D}$
of $\xi$, where $D$ is a non-singular odd theta-characteristic.
Then the identity (\ref{quasidetid}) can be rewritten as
$$
F(\chi(x_1-y_1),x_0,y_0)\frac{\th_{D}(x_0-x_1)\th_{D}(y_0-y_1)}
{\th_{D}(x_0-y_1)\th_{D}(y_0-x_1)}
=F(\chi,x_0,y_0)-F(\chi,x_1,y_0)\circ
F(\chi,x_1,y_1)^{-1}\circ F(\chi,x_0,y_1).
$$
Note that the fraction in the left-hand side is equal to
$\frac{E(x_0,x_1)E(y_0,y_1)}{E(x_0,y_1)E(y_0,x_1)}$.
Similarly, the identity of Theorem \ref{quasidetthm} implies that
$$|F(\chi,x_j,y_i)_{0\le i,j\le n}|_{00}=
F(\chi(\sum_{i=1}^n x_i-\sum_{i=1}^n y_i),x_0,y_0)\cdot
\prod_{i=1}^n\frac{E(x_0,x_i)E(y_0,y_i)}{E(x_0,y_i)E(y_0,x_i)}.$$
This is precisely the identity appearing in the proof of
Lemma 2 in \cite{F3} (instead of quasideterminants, Fay uses the inverse
of a matrix with non-commuting entries).

\section{Tangents to the canonical embedding}

In this section we will relate the identity of Corollary
\ref{cor2} to tangent lines to a canonically embedded curve.

Let us rewrite this Corollary slightly.
Pick generic elements $\eta_i\in H^0(C,\om_C(-E_i-t))$ for $i=1,2$.
Then we can easily calculate all the Massey products appearing in
Corollary \ref{cor2}. For example, since
$\eta_2/\eta$ is a rational function with poles at $E_1+x+y+z$,
we get
$$m_3(\OO_C(E_1+y+z),x,z)=\frac{\eta_2(z)\eta'(x)}{\eta_2(x)\eta'(z)},$$
etc. Thus, the identity of Corollary \ref{cor2} becomes
$$\frac{\eta_1(x)\eta_2(x)}{\eta'(x)}+\frac{\eta_1(y)\eta_2(y)}{\eta'(y)}+
\frac{\eta_1(z)\eta_2(z)}{\eta'(z)}=0$$
which is nothing else but the residue theorem applied to
$\eta_1\eta_2/\eta$. This is not surprizing since the
proof of Theorem \ref{idthm} was based on the residue theorem.
Of course, similar relation holds for arbitrary
global $1$-forms $\eta_1$ and $\eta_2$. However, one may expect that
in the case when $\eta_1\eta_2$ vanishes at most of the zeroes of
$\eta$, such a relation is more valuable. 

Now let us assume that $C$ is non-hyperelliptic and consider
the canonical embedding $C\hookrightarrow{\Bbb P}(H^0(C,\om_C)^*)$.
Then $\eta$, $\eta_1$ and $\eta_2$ are just equations of hyperplanes.
To interpret the derivative $\eta'(x)$ we note that
there is a canonical isomorphism $\om_C\simeq\OO(1)|_C$. Therefore,
a lifting of a point $x\in C\subset{\Bbb P}(H^0(C,\om_C)^*)$
to a (non-zero) vector
$\wt{x}\in H^0(C,\om_C)^*$ determines canonically a tangent
vector $t_{\wt{x}}$ to $C$ at $x$, which depends linearly on $\wt{x}$.
Recall that the tangent space to a point $L$ of a projective
space $\Proj(V)$ is isomorphic to $\Hom(L,V/L)$. Thus, a tangent
vector $t_{\wt{x}}$ determines an element
$v_{\wt{x}}:=t_{\wt{x}}(\wt{x})\in H^0(C,\om_C)^*/x$ depending
quadratically on $\wt{x}$. If $Q$ is a quadratic form on 
the space $H^0(C,\om_C)^*$ and $l$ is a linear form on this space
such that $l(x)=0$, then the ratio
$Q(\wt{x})/l(v_{\wt{x}})$ depends only on $x$. By abuse of notation
we denote this ratio by $Q(x)/l(v_x)$.

\begin{prop}\label{canprop} Let $C$ be a non-hyperelliptic curve,
$D$ be the intersection of $C$
with a hyperplane $l=0$ in the canonical embedding.
Assume that $D$ consists of $2g-2$ distinct points.
Then for any quadratic form $Q$ on $H^0(C,\om_C)^*$ one has
$$\sum_{x\in D}\frac{Q(x)}{l(v_x)}=0.$$
\end{prop}

\Pf . It suffices to prove this in the case when $Q$ is a product
of linear forms. Then it follows immediately from the residue theorem
as we have seen above.
\ed

Using the above proposition we can compute the tangent line
to $C$ at a given point $x$ in terms of some ratios associated
with hyperplane sections of $C$ containing $x$.
Note that the tangent line to $C$ at $x$
defines a point in the $(g-2)$-dimensional projective space which is
the target of the linear projection from $x$. We start by computing the
image of this tangent line under the
linear projection with the center at the
line spanned by $x$ and another point $y\in C$.
Our calculation is based on the following result.

\begin{cor}\label{ratiocor}
Let $x\neq y$ be a pair of generic points on a canonically
embedded curve $C$, $\wt{x},\wt{y}\in H^0(C,\om_C)^*$
be their liftings. Let $H$ be a generic hyperplane passing through
$x,y$. Write the intersection $C\cap H$ in the form
$x+y+D_1+D_2$ where $D_1$ and $D_2$ are divisors of degree $g-2$.
Let $l_{D_i}=0$ for $i=1,2$ be an equation of the unique hyperplane
in $H$ passing through $D_i$. Then the ratio
$$r(\wt{x},\wt{y},H)=\frac{l_{D_1}(\wt{y})l_{D_2}(\wt{y})}
{l_{D_1}(\wt{x})l_{D_2}(\wt{x})}$$
does not depend on a choice of $D_1$ and $D_2$.
More precisely,
$$r(\wt{x},\wt{y},H)=-\frac{l(v_{\wt{y}})}{l(v_{\wt{x}})}$$
where $l=0$ is an equation of $H$.
\end{cor}

\Pf . Take $Q=l_{D_1}l_{D_2}$ in Proposition \ref{canprop}.  
\ed

Below we denote by $(t_1:\ldots:t_n)$
homogeneous coordinates of a point in $\Proj^{n-1}$.
Also for a linear form $l$ we set
$r(\wt{x},\wt{y},l):=r(\wt{x},\wt{y},(l=0))$.

\begin{cor} Let $x,y\in C$ be points as in the Corollary \ref{ratiocor},
$\wt{x},\wt{y}\in H^0(C,\om_C)^*$ be their liftings.
Pick $g-2$ generic linear forms $l_1,\ldots,l_{g-2}$
vanishing at $x,y$. Then one has
\begin{align*}
&(l_1(v_{\wt{x}}):\ldots:l_{g-2}(v_{\wt{x}}))=\\
&(1:\frac{b_2-b_1}{a_2-b_2}:\frac{(a_2-b_1)(b_3-b_2)}{(a_2-b_2)(a_3-b_3)}:
\ldots:\frac{(a_2-b_1)\ldots(a_{g-3}-b_{g-4})(b_{g-2}-b_{g-3})}
            {(a_2-b_2)\ldots(a_{g-3}-b_{g-3})(a_{g-2}-b_{g-2})}),
\end{align*}
where $a_i=r(\wt{x},\wt{y},l_i)$, $b_i=r(\wt{x},\wt{y},l_1+\ldots+l_i)$.
\end{cor}

\Pf . Set $u_i=l_i(v_{\wt{x}})$, $v_i=l_i(v_{\wt{y}})$,
where $i=1,\ldots,g-2$.
According to Corollary \ref{ratiocor}, we have
$$\frac{v_i}{u_i}=-a_i, \ i=1,\ldots,g-2;$$
$$\frac{v_1+\ldots+v_i}{u_1+\ldots+u_i}=-b_i,\ i=1,\ldots,g-2.$$
It follows that
$$\frac{u_1+\ldots+u_i}{u_{i+1}}=\frac{a_{i+1}-b_{i+1}}{b_{i+1}-b_i},
\ i=1,\ldots,g-3.$$
From this relation we easily deduce that
$$\frac{u_{i+1}}{u_i}=\frac{(a_i-b_{i-1})(b_{i+1}-b_i)}
{(b_i-b_{i-1})(a_{i+1}-b_{i+1})}$$
for $i\ge 2$ while 
$$\frac{u_2}{u_1}=\frac{b_2-b_1}{a_2-b_2}.$$
\ed

Applying the above Corollary to the pairs $(x,p)$ and $(x,q)$, where
$p$ and $q$ are generic points on $C$, 
we obtain the following formula for the tangent line to $C$ at
$x$.

\begin{cor} Let $x$ be a point on $C$, $\wt{x}\in H^0(C,\om_C)^*$ be
its lifting. Let $p,q\in C$ be a pair of generic points, $\wt{p}$,
$\wt{q}$ be their liftings. Pick
$g-3$ generic linear forms $l_1,\ldots,l_{g-3}$ vanishing at $x$, $p$
and $q$. Then pick a generic linear form $l_{g-2}$ vanishing at $x$ and $p$
and a generic linear form $l_0$ vanishing at $x$ and $q$.
Then 
\begin{align*}
&(l_0(v_{\wt{x}}):l_1(v_{\wt{x}}):\ldots:l_{g-2}(v_{\wt{x}}))=\\
&(\frac{c_1-d}{d-c_0}:1:
\frac{b_2-b_1}{a_2-b_2}:\frac{(a_2-b_1)(b_3-b_2)}{(a_2-b_2)(a_3-b_3)}:
\ldots:\frac{(a_2-b_1)\ldots(a_{g-3}-b_{g-4})(b_{g-2}-b_{g-3})}
            {(a_2-b_2)\ldots(a_{g-3}-b_{g-3})(a_{g-2}-b_{g-2})}),
\end{align*}
where $a_i=r(\wt{x},\wt{p},l_i)$, $b_i=r(\wt{x},\wt{p},l_1+\ldots+l_i)$,
$c_i=r(\wt{x},\wt{q},l_i)$, $d=r(\wt{x},\wt{q},l_0+l_1)$.
\end{cor}

\vspace{3mm}

{\sc Department of Mathematics and Statistics, 
Boston University, Boston MA 02215}

{\it E-mail address:} apolish@@math.bu.edu

\end{document}